\documentclass[12pt,letterpaper,twoside]{article}

\usepackage[centertags]{amsmath}
\usepackage{amsmath,amsfonts,amssymb,amsthm,newlfont}
\font\Goth=yinitas scaled \magstep0 
\newcommand{\Gth}[1]{\lower2mm\hbox{\Goth #1}}

\usepackage[latin1]{inputenc}
\def\l1{{\lambda}_1}

\newcommand{\f}{\frac}

\def\x1{{\xi }_{xx}}
\def\x2{{\xi }_{yy}}
\def\x3{{\xi }_{xy}}

\def\e1{{\eta }_{xx}}
\def\e2{{\eta }_{yy}}
\def\e3{{\eta }_{xy}}

\newtheorem{theorem}{Theorem}

\newcommand{\beqn}{\begin{eqnarray*}}
\newcommand{\eeqn}{\end{eqnarray*}}
\newcommand{\beqnn}{\begin{eqnarray}}
\newcommand{\eeqnn}{\end{eqnarray}}
\newcommand{\p}{\partial}

\newcommand{\bb}{\begin{equation}}
\newcommand{\ee}{\end{equation}}
\newcommand{\ba}{\begin{array}}
\newcommand{\ea}{\end{array}}
\newcommand{\R}{\mathbb{R}}

\newcommand{\Hi}{\mathbb{H}}
\newcommand{\lh}{\Delta_{\Hi}}
\newcommand{\X}{\tilde{X}}
\newcommand{\Y}{\tilde{Y}}

\begin{document}
\pagenumbering{arabic}
\title{\huge \bf Conservations Laws for Critical Kohn-Laplace Equations on the Heisenberg group}
\author{\rm \large Yuri Bozhkov and Igor Leite Freire\\ \\
\it Instituto de Matem\'atica,
Estat\'\i stica e \\ \it Computa\c c\~ao Cient\'\i fica - IMECC \\
\it Universidade Estadual de Campinas - UNICAMP \\ \it C.P.
$6065$, $13083$-$970$ - Campinas - SP, Brasil
\\ \rm E-mail: bozhkov@ime.unicamp.br \\ \ \ \ \ \
igor@ime.unicamp.br }
\date{\ }
\maketitle \vspace{1cm}
\begin{abstract}
Using the complete group classification of semilinear differential
equations on the three-dimensional Heisenberg group $\Hi$, carried
out in a preceding work, we establish the conservation laws for
the critical Kohn-Laplace equations via the Noether's Theorem.
\end{abstract}

\vskip 1cm
\begin{center}
{2000 AMS Mathematics Classification numbers:\vspace{0.2cm}\\
35H10, 58J70\vspace{0.2cm} }
\end{center}

\pagenumbering{arabic}
\parindent 30pt
\baselineskip 16pt

\section{Introduction}

In a previous work \cite{yi1}, we obtained the complete group
classification of the following semilinear equation on the
three-dimensional Heisenberg group $\Hi$:
\begin{equation}\label{kl}
\lh u+f(u)=0,
\end{equation}
(Here $\lh$ is the Kohn-Laplace operator.). Further, we showed in
\cite{yi2} that all Lie point symmetries of (\ref{kl}) in the
critical Stein-Sobolev case $f(u)=u^{3}$ are variational or
divergence symmetries.

The purpose of this note is to establish the corresponding
conservation laws via the Noether's Theorem (\cite{bk},\cite{ol}).
As it is well known, the latter provides an algorithimic procedure
for construction of conservation laws. Namely,
let\begin{equation*} X=\xi^{i}\f{\p}{\p x^{i}}+\eta\f{\p}{\p u}
\end{equation*}
be the generator of an infinitesimal transformation admitted by
certain Euler-Lagrange equation $E(\cal{L})=0$ of order $2k$,
whose Lagrangian is denoted by $\cal{L}$. If X is a divergence
symmetry of $E(\cal{L} )=0$, that is, if there exists a vector
valued function $\varphi=(\varphi^{i})$  such that
\begin{equation}\label{divs}
X^{(k)}\cal{L}+\cal{L}D_{i}\xi^{i}=D_{i}\varphi^{i},
\end{equation}
then the Noether's Theorem states that the following conservation
law holds
\begin{equation}\label{cons}
D_{i}(\xi^{i}\cal{L}+W^{i}[u,\eta-\xi^{j}u_{j}]-\varphi^{i})=0
\end{equation}
for all solutions $u$ of $E(\cal{L})=0$. Above we have used the
same notations and conventions as in \cite{bk}. (For the
definition of $W^{i}$ see \cite{bk}, pp. 254-255.) Therefore, as
pointed out in \cite{bk}, to apply this theorem one must

(i) find all transformations admitted by $E(\cal{L})=0$ and

(ii) check which infinitesimal generators $X$ satisfy the
condition (\ref{divs}).

Hence it is clear that the major difficulty in applying the
Noether's Theorem is that usually there is no explicit formula for
the potential $\varphi$.

As it was already mentioned in the begining, the first step (i) is
done in \cite{yi1}  and the second (ii) - in the work \cite{yi2}.
Moreover, we observe that in \cite{yi2} we found explicity the
potentials $\varphi$ associated to the corresponding divergence
symmetries of the equation
\begin{equation}\label{crit}
\lh u+u^{3}=0.
\end{equation}

Thus we have at our disposal all ingredients which will enable us
to apply \textit{directly} the Noether's Theorem by a
straightforward calculation.

In this paper we are interested in the critical Kohn-Laplace
equation (\ref{crit}) since it possesses the widest symmetry group
among the {\it nonlinear} equations of form (\ref{kl}). See
\cite{yi1}. The Noether symmetries and the corresponding
conservation laws in the {\it linear} cases $f(u)=0$ and $f(u)=u$
will be treated elsewhere.

The next step in this research is to construct nonlocal symmetries
and the corresponding to them nonlocal conservation laws for the
solutions of critical semilinear Kohn-Laplace equations on the
Heisenberg Group using the recent methods devised and developed by
George Bluman et al. This problem will be treated elsewhere. Here
we merely point out that the obtained (local) conservation laws
will be used for that purpose.

The paper is organized as follows. In the section 2 we present
briefly some of the main aspects of Heisenberg groups as well as
parts of the results, obtained in (\cite{yi1},\cite{yi2}), which
will be used later. The conservation laws are stated in section 3,
in the form of Theorem 1.

\section{The Noether symmetries of critical Kohn-Laplace equa-tions}

To begin with, we recall some facts concerning the Heisenberg
group $\Hi$. Let $\phi:\R^{3}\times\R^{3}\rightarrow\R^{3}$,
defined by $\phi((x,y,t),(x_{0},y_{0},t_{0})):=(x+x_{0}, y+y_{0},
t+t_{0}+2(xy_{0}-yx_{0})),$ be the composition law of $\Hi$
determining its Lie group structure. The following vector fields
\begin{equation}\label{fields}
\begin{array}{l c c c l}
X & = &\displaystyle{ \f{d}{ds}\phi((x,y,t),(s,0,0))|_{s=0}} & = &\displaystyle{ \f{\p}{\p x}+2y\f{\p}{\p t}},\\
\\
Y & = &\displaystyle{ \f{d}{ds}\phi((x,y,t),(0,s,0))|_{s=0}} & = & \displaystyle{\f{\p}{\p y}-2x\f{\p}{\p t}},\\
\\
T & =& \displaystyle{ \f{d}{ds}\phi((x,y,t),(0,0,s))|_{s=0}} & = &
\displaystyle{\f{\p}{\p t}}
\end{array}
\end{equation}
form a basis of left invariant vectors fields on $\Hi$. The
Riemannian metric $ds^{2}=dx^{2}+dy^{2}+(2ydx-2xdy+dt)^{2}$ is a
left invariant metric and the Lie algebra generated by
\begin{equation}\label{isovec} T=\frac{\p }{\p t},\;\;\; R= y\frac{\p }{\p
x}-x \frac{\p }{\p y},\;\;\; \X=\frac{\p}{\p x}-2y \frac{\p}{\p
t}, \;\;\;\Y=\frac{\p}{\p y}+2x \frac{\p}{\p t} \end{equation} is
the Lie algebra of the infinitesimal isometries of $\Hi$.

The Kohn-Laplace operator is defined by $\lh: =X^2+Y^2$, where $X$
and $Y$ are defined in (\ref{fields}). For $u=u(x,y,t): {\R
}^3\rightarrow \R $, we have $ \lh
u=u_{xx}+u_{yy}+4(x^{2}+y^{2})u_{tt}+4yu_{xt}-4xu_{yt} . $

We point out that the Kohn-Laplace operator $\lh $ is not a
(strongly) elliptic operator. Nevertheless it was shown in
(\cite{yi1},\cite{yi2}) that the Lie symmetry theory can be
successfully applied to such a {\it subelliptic} operator.

The equation (\ref{crit}) arises from the following Lagrangian
\begin{equation}\label{lag}
\cal{L}=\frac{1}{2}u_{x}^{2}+\frac{1}{2}u_{y}^{2}+2(x^{2}+y^{2})u_{t}^{2}+2yu_{x}u_{t}-2xu_{y}u_{t}-\f{u^{4}}{4}.
\end{equation}

By the group classification \cite{yi1}, the symmetry algebra of
(\ref{crit}) is generated by (\ref{isovec}) and the following
vectors fields:
\begin{equation*}
\begin{array}{l}
\displaystyle{Z= x\frac{\p}{\p x}+y\frac{\p}{\p y}+2t
\frac{\p}{\p t}-u\frac{\p}{\p u}}, \\
\\
 \displaystyle{V_1 = (xt-x^{2}y-y^{3})\frac{\p }{\p x} +
(yt+x^{3}+xy^{2})\frac{\p }{\p y} +
(t^{2}-(x^{2}+y^{2})^{2})\frac{\p }{\p t}-
t u \frac{\p }{\p u}},\\
\\
\displaystyle{V_2  = (t-4xy)\frac{\p }{\p x} +
(3x^{2}-y^{2})\frac{\p }{\p y} - (2yt+2x^{3}+2xy^{2})\frac{\p }{\p
t} + 2 y u \frac{\p }{\p u}},
\\
\\
\displaystyle{V_3 = (x^{2}-3y^{2})\frac{\p }{\p x} +
(t+4xy)\frac{\p }{\p y} + (2xt-2x^{2}y-2y^{3})\frac{\p }{\p t} -
2x u \frac{\p }{\p u}}.
\end{array}
\end{equation*}

In this case we have the following comutation table, not presented
in \cite{yi1}.
\begin{table}[htbp]
\begin{center}
\begin{tabular}{|c|c|c|c|c|c|c|c|c|}\hline
        & T &  R     & $\X$        & $\Y$         &  $V_{1}$   & $V_{2}$ & $V_{3}$ & $Z$       \\\hline
T       & 0 &  0     &     0       &   0          &  $Z$ & $\X$
& $\Y$    & 2T            \\\hline R       & 0 &  0     & $\Y$
& -$\X$        &  0         & $V_{3}$ & -$V_{2}$& 0
\\\hline $\X$    & 0 &-$\Y$   &     0       &    4T        &
$V_{2}$   & -6R  &$2Z$ & $\X$         \\\hline $\Y $   & 0 &$\X$
&    - 4T      &     0        & $V_{3}$ &$-2Z$&-6R      & $\Y$
\\\hline $V_{1}$ & $-Z$ &  0  & -$V_{2}$    &     $-V_{3}$  &  0
&0        & 0       &-2$V_{1}$\\ \hline $V_{2}$ & -$\X$      &
-$V_{3}$   &    6R &  2Z   &  0     &  0      &4$V_{1}$ &
-$V_{2}$\\ \hline $V_{3}$ & $-\Y$&-$V_{2}$ & -2Z & $6R$    &  0
&-4$V_{1}$&     0   &   -$V_{3}$\\ \hline $Z$ & -2T        &  0
&    -$\X$ &  -$\Y$&  2$V_{1}$  & $V_{2}$ & $V_{3}$ &   0
\\ \hline
\end{tabular}
\caption{\small{Table of Lie brackets of equation (\ref{crit})}}
\end{center}
\end{table}

In \cite{yi2} we showed that (\ref{isovec}) and $Z$ are
variational symmetries ($\varphi=0$ in (\ref{divs})) and $V_{1}$,
$V_{2}$ and $V_{3}$ are divergence symmetries of (\ref{crit}).
Hence all Lie point symmetries of the critical Kohn-Laplace
equation (\ref{crit}) are Noether symmetries (variational or
divergence symmetries). Moreover, in \cite{yi2} we found
explicitly the potentials $\varphi$ in the conservation laws
implied by the Noether's Theorem.

In the next section we state the main result of this paper.

\newpage
\section{The Conservation Laws}

\begin{theorem}
The conservations laws of the Noether symmetries are:
\begin{enumerate}

\item For the symmetry $T$, the conservation law is $Div(\tau)=0$,
where $\tau=(\tau_{1},\tau_{2},\tau_{3})$ and
\begin{equation*}
\begin{array}{l}
\displaystyle{\tau_{1} = -2yu_{t}^{2}-u_{x}u_{t}},
\\
\\
\displaystyle{\tau_{2} = 2xu_{t}^{2}-u_{y}u_{t}},
\\
\\
\displaystyle{\tau_{3} =
\f{1}{2}u_{x}^{2}+\f{1}{2}u_{y}^{2}-2(x^{2}+y^{2})u_{t}^{2}-\frac{1}{4}u^{4}}.
\end{array}
\end{equation*}

\item For the symmetry $R$, the conservation law is
$Div(\sigma)=0$, where $\sigma=(\sigma_{1},\sigma_{2},\sigma_{3})$
and
\begin{equation*}
\begin{array}{l}
\displaystyle{\sigma_{1} =
-\f{1}{2}yu_{x}^{2}+\f{1}{2}yu_{y}^{2}+2y(x^{2}+y^{2})u_{t}^{2}+xu_{x}u_{y}-\frac{1}{4}yu^{4}},
\\
\\
\displaystyle{\sigma_{2} =
-\f{1}{2}xu_{x}^{2}-\f{1}{2}xu_{y}^{2}-2x(x^{2}+y^{2})u_{t}^{2}-yu_{x}u_{y}+\frac{1}{4}xu^{4}},
\\
\\
\displaystyle{\sigma_{3} =
-2y^{2}u_{x}^{2}-2x^{2}u_{y}^{2}+4xyu_{x}u_{y}-4y(x^{2}+y^{2})u_{x}u_{t}+4x(x^{2}+y^{2})u_{y}u_{t}}.
\end{array}
\end{equation*}

\item For the symmetry $\tilde{X}$, the conservation law is
$Div(\chi)=0$, where $\chi=(\chi_{1},\chi_{2},\chi_{3})$ and
\begin{equation*}
\begin{array}{l}
\displaystyle{\chi_{1} =
-\f{1}{2}u_{x}^{2}+\f{1}{2}u_{y}^{2}+2(x^{2}+3y^{2})u_{t}^{2}+2yu_{x}u_{t}-2xu_{y}u_{t}-\frac{1}{4}u^{4}},
\\
\\
\displaystyle{\chi_{2} =
-4xyu_{t}^{2}-u_{x}u_{y}+2xu_{x}u_{t}+2yu_{y}u_{t}},
\\
\\
\displaystyle{\chi_{3}
=-3yu_{x}^{2}-yu_{y}^{2}+4y(x^{2}+y^{2})u_{t}^{2}+2xu_{x}u_{y}-4(x^{2}+y^{2})u_{x}u_{t}+\frac{1}{2}yu^{4}}.
\end{array}
\end{equation*}

\item For the symmetry $\tilde{Y}$, the conservation law is
$Div(\upsilon)=0$, where
$\upsilon=(\upsilon_{1},\upsilon_{2},\upsilon_{3})$ and
\begin{equation*}
\begin{array}{l}
\displaystyle{\upsilon_{1} =
-4xyu_{t}^{2}-u_{x}u_{y}-2xu_{x}u_{t}-2yu_{y}u_{t}},
\\
\\
\displaystyle{\frac{1}{2}u_{x}^{2}-\frac{1}{2}u_{y}^{2}+2(3x^{2}+y^{2})u_{t}^{2}+2yu_{x}u_{t}-2xu_{y}u_{t}-\frac{1}{4}u^{4}},
\\
\\
\displaystyle{\upsilon_{3} =
xu_{x}^{2}+3xu_{y}^{2}-4x(x^{2}+y^{2})u_{t}^{2}-2yu_{x}u_{y}-4(x^{2}+y^{2})u_{y}u_{t}-\frac{1}{2}xu^{4}}.
\end{array}
\end{equation*}

\item For the symmetry $Z$, the conservation law is
$Div(\zeta)=0$, where $\zeta=(\zeta_{1},\zeta_{2},\zeta_{3})$ and
\begin{equation*}
\begin{array}{l}
\begin{array}{l c l}
\zeta_{1} &= &\displaystyle{-\f{1}{2}xu_{x}^{2}+\f{1}{2}xu_{y}^{2}+2(x^{3}-2ty+xy^{2})u_{t}^{2}-y u_{x}u_{y}-2t u_{x}u_{t}}\\
\\
&
&\displaystyle{-2(x^{2}+y^{2})u_{y}u_{t}-uu_{x}-2yuu_{t}-\frac{1}{4}x
u^{4}},
\end{array}
\\
\\
\begin{array}{l c l}
\zeta_{2}& = & \displaystyle{\frac{1}{2}y u_{x}^{2}-\frac{1}{2}y u_{y}^{2}+2(2tx+x^{2}y+y^{3})u_{t}^{2}-x u_{x}u_{y}+2(x^{2}+y^{2})u_{x}u_{t}}\\
\\
& & \displaystyle{-2t u_{y}u_{t}-uu_{y}+2xuu_{t}-\frac{1}{4}y
u^{4}},
\end{array}
\\
\\
\begin{array}{l c l}
\zeta_{3} & = & \displaystyle{(t-2xy)u_{x}^{2}+(t+2xy)u_{y}^{2}-4t(x^{2}+y^{2})u_{t}^{2}+2(x^{2}-y^{2})u_{x}u_{y}-4x(x^{2}+y^{2})u_{x}u_{t}}\\
\\
&
&\displaystyle{-4y(x^{2}+y^{2})u_{y}u_{t}+2xuu_{y}-2yuu_{x}-4(x^{2}+y^{2})uu_{t}-\frac{1}{2}t
u^{4}}.
\end{array}
\end{array}
\end{equation*}

\item For the symmetry $V_{1}$, the conservation law is
$Div(A)=0$, where $A=(A_{1},A_{2},A_{3})$ and
\begin{equation*}
\begin{array}{l}
\begin{array}{l c l}
A_{1}& = &\displaystyle{-\frac{1}{2}(tx-x^{2}y-y^{3})u_{x}^{2}+\frac{1}{2}(tx-x^{2}y-y^{3})u_{y}^{2}+2t(x^{3}+xy^{2}-ty)u_{t}^{2}}\\
\\
& &\displaystyle{-(x^{3}+xy^{2}+ty)u_{x}u_{y}-[t^{2}-(x^{2}+y^{2})^{2}]u_{x}u_{t}-2t(x^{2}+y^{2})u_{y}u_{t}}\\
\\
&
&\displaystyle{-tuu_{x}-2tyuu_{t}+yu^{2}-\frac{1}{4}(tx-x^{2}y-y^{3})u^{4}},
\end{array}
\\
\\
\begin{array}{l c l}
A_{2}& = & \displaystyle{\frac{1}{2}(x^{3}+ty+xy^{2})u_{x}^{2}-\frac{1}{2}(x^{3}+ty+xy^{2})u_{y}^{2}+2t(x^{2}y+y^{3}+tx)u_{t}^{2}}\\
\\
& & \displaystyle{-(tx-x^{2}y-y^{3})u_{x}u_{y}+2t(x^{2}+y^{2})u_{x}u_{t}-[t^{2}-(x^{2}+y^{2})^{2}]u_{y}u_{t}}\\
\\
& &
\displaystyle{-tuu_{y}+2txuu_{t}-xu^{2}-\frac{1}{4}(x^{3}+ty+xy^{2})u^{4}},
\end{array}
\\
\\
\begin{array}{l c l}
A_{3} & = & \displaystyle{+\frac{1}{2}(t^{2}-x^{4}-4txy+2x^{2}y^{2}+3y^{4})u_{x}^{2}+\frac{1}{2}(t^{2}+3x^{4}+4txy+2x^{2}y^{2}-y^{4})u_{y}^{2}}\\
\\
& &\displaystyle{-2(x^{2}+y^{2})[t^{2}-(x^{2}+y^{2})^{2}]u_{t}^{2}+2[t(x^{2}-y^{2})-2xy(x^{2}+y^{2})]u_{x}u_{y}}\\
\\
& &\displaystyle{-4(x^{2}+y^{2})(tx-x^{2}y-y^{3})u_{x}u_{t}-4(x^{2}+y^{2})(x^{3}+ty+xy^{2})u_{y}u_{t}}\\
\\
&
&\displaystyle{-2tyuu_{x}+2txuu_{y}-4t(x^{2}+y^{2})uu_{t}+2(x^{2}+y^{2})u^{2}-\frac{1}{4}[t^{2}-(x^{2}+y^{2})^{2}]u^{4}}.\end{array}
\end{array}
\end{equation*}
\newpage

\item For the symmetry $V_{2}$, the conservation law is
$Div(B)=0$, where $B=(B_{1},B_{2},B_{3})$ and
\begin{equation*}
\begin{array}{l}
\begin{array}{l c l}
B_{1}& = &
\displaystyle{-\frac{1}{2}(t-4xy)u_{x}^{2}+\frac{1}{2}(t-4xy)u_{y}^{2}+[2t(x^{2}+3y^{2})-4xy(x^{2}+y^{2})]u_{t}^{2}}
\\
\\
& &+\displaystyle{-(3x^{2}-y^{2})u_{x}u_{y}+2(x^{3}+ty+xy^{2})u_{x}u_{t}-2(tx-x^{2}y-y^{3})u_{y}u_{t}}\\
\\
& & \displaystyle{+2yuu_{x}+4y^{2}uu_{t}-\frac{1}{4}(t-4xy)u^{4}},
\end{array}
\\
\\
\begin{array}{l c l}
B_{2} & = & \displaystyle{\frac{1}{2}(3x^{2}-y^{2})u_{x}^{2}-\frac{1}{2}(3x^{2}-y^{2})u_{y}^{2}+2(x^{4}-2txy-y^{4})u_{t}^{2}-(t-4xy)u_{x}u_{y}}\\
\\
& &\displaystyle{+2(tx-x^{2}y-y^{3})u_{x}u_{t}+2(x^{3}+ty+xy^{2})u_{y}u_{t}+2yuu_{y}-4xyuu_{t}-u^{2}}\\
\\
& &\displaystyle{-\frac{1}{4}(3x^{2}-y^{2})u^{4}},
\end{array}
\\
\\
\begin{array}{l c l}
B_{3} & = & \displaystyle{(7xy^{2}-x^{3}-3ty)u_{x}^{2}+(5x^{3}-3xy^{2}-ty)u_{y}^{2}+4(x^{2}+y^{2})(x^{3}+ty+xy^{2})u_{t}^{2}}\\
\\
& &\displaystyle{+2(tx-7x^{2}y+y^{3})u_{x}u_{y}-4(t-4xy)(x^{2}+y^{2})u_{x}u_{t}-4(3x^{4}+2x^{2}y^{2}-y^{4})u_{y}u_{t}}\\
\\
& &\displaystyle{ +2xu^{2}+4y^{2}uu_{x}-4xyuu_{y}+8y(x^{2}+y^{2})uu_{t}+\frac{1}{2}(x^{3}+ty+xy^{2})u^{4}}\\
\end{array}
\end{array}
\end{equation*}

\item For the symmetry $V_{3}$, the conservation law is
$Div(C)=0$, where $C=(C_{1},C_{2},C_{3})$ and
\begin{equation*}\begin{array}{l}
\begin{array}{l c l}
C_{1}& =& \displaystyle{\frac{1}{2}(x^{2}y-tx+y^{3})u_{x}^{2}+\frac{1}{2}(tx-x^{2}y-y^{3})u_{y}^{2}+2t(x^{3}-ty+xy^{2})u_{t}^{2}}\\
\\
& &\displaystyle{-(x^{3}+ty+xy^{2})u_{x}u_{y}-[t^{2}-(x^{2}+y^{2})^{2}]u_{x}u_{t}-2t(x^{2}+y^{2})u_{y}u_{t}}\\
\\
&
&\displaystyle{-tuu_{x}-2tyuu_{t}-\frac{1}{4}(tx-x^{2}y-y^{3})u^{4}},
\end{array}
\\
\\\begin{array}{l c l}
C_{2}&=&\displaystyle{\frac{1}{2}(x^{3}+ty+xy^{2})u_{x}^{2}-\frac{1}{2}(x^{3}+ty+xy^{2})u_{y}^{2}+2t(tx+x^{2}y+y^{3})u_{t}^{2}}\\
\\
& & \displaystyle{-(tx-x^{2}y-y^{3})u_{x}u_{y}+2t(x^{2}+y^{2})u_{x}u_{t}-[t^{2}-(x^{2}+y^{2})^{2}]u_{y}u_{t}}\\
\\
&
&\displaystyle{-u^{2}-tuu_{y}+2txuu_{t}-\frac{1}{4}(x^{3}+ty+xy^{2})u^{4}},
\end{array}
\end{array}
\end{equation*}

\begin{equation*}
\begin{array}{l}
\begin{array}{l c l}
C_{3}&=&\displaystyle{\frac{1}{2}(t^{2}-x^{4}-4txy+2x^{2}y^{2}+3y^{4})u_{x}^{2}+\frac{1}{2}(t^{2}+3x^{4}+4txy+2x^{2}y^{2}-y^{4})u_{y}^{2}}\\
\\
& &\displaystyle{-2(x^{2}+y^{2})[t^{2}-(x^{2}+y^{2})^{2}]u_{t}^{2}+2[t(x^{2}-y^{2})-2xy(x^{2}+y^{2})]u_{x}u_{y}}\\
\\
& &\displaystyle{+4(x^{2}+y^{2})(x^{2}y-tx+y^{3})u_{x}u_{t}-4(x^{2}+y^{2})(x^{3}+ty+xy^{2})u_{y}u_{t}}\\
\\
&
&\displaystyle{+2txuu_{y}-2tyuu_{x}-4t(x^{2}+y^{2})uu_{t}+2yu^{2}-\frac{1}{4}[t^{2}-(x^{2}+y^{2})^{2}]u^{4}}.
\end{array}
\end{array}
\end{equation*}
\end{enumerate}
\end{theorem}

\begin{proof} First, we observe that the potentials for the
symmetries $T,\;R,\;\X,\;\Y$ and Z of (\ref{crit}) are 0, that is,
these symmetries are variational \cite{yi2}. Further, the
potentials $\varphi$ of $V_{1}$, $V_{2}$, $V_{3}$ are
$(-yu^2,xu^2,-2(x^2+y^2)u^2), (0,u^2,-2x u^2), (-u^2,0,-2y u^2)$
respectively. See \cite{yi2}.

Now, with these at hand, as mentioned in the introduction, the
proof is by a tedious straightforward calculation, which we shall
not present here for obvious reasons. However, a computer assisted
proof can be obtained by two simple {\it Mathematica} programs.
The first one calculates the components of the conservation laws,
which appear in the equation (\ref{cons}). The second program
verifies the conservation laws using the Noether Identity
\cite{ib}. Both Mathematica notebooks can be obtained form the
authors upon request.\end{proof}

\begin{center}\textbf{Acknowledgements}\end{center}

 We wish to thank Lab. EPIFISMA (Proj. FAPESP) for having given us
 the opportunity to use excellent computer facilities. Y. Bozhkov would also like to
 thank CNPq, Brasil, for partial financial support. I. L. Freire is grateful to
 UNICAMP for financial support.

 \end{document}